\documentclass{amsart}

\usepackage{amssymb,amsmath,amsthm,color}

\usepackage[bookmarksopen=true,bookmarksnumbered=true, bookmarksopenlevel=2, colorlinks=true,linkcolor=mblue,
citecolor=mblue,urlcolor=mblue, pdfstartview=FitH]{hyperref}

\definecolor{mblue}{rgb}{0,0,.8}

\usepackage[all]{xy}

\newcommand{\Z}{\mathbb Z}
\newcommand{\Q}{\mathbb Q}
\newcommand{\F}{\mathbb F}

\newcommand{\C}{\mathbb C}
\newcommand{\p}{\mathfrak p}
\newcommand{\e}{\epsilon}
\newcommand{\W}{\mathfrak W}

\newtheorem{thm}{Theorem}
\newtheorem{lem}{Lemma}

\newtheorem{prop}{Proposition}
\newtheorem{cor}{Corollary}

\DeclareMathOperator{\PGL}{PGL}
\DeclareMathOperator{\PSL}{PSL}
\DeclareMathOperator{\GL}{GL}
\DeclareMathOperator{\Hc}{H}
\DeclareMathOperator{\Ima}{Im}
\DeclareMathOperator{\Gal}{Gal}
\DeclareMathOperator{\Ind}{Ind}
\DeclareMathOperator{\Proj}{Proj}

\def\dash---{\thinspace---\hskip.16667em\relax}

\begin{document}

\title{On modular mod $\ell$ Galois representations with exceptional images.}
\author[Ian Kiming, Helena A. Verrill]{Ian Kiming, Helena A. Verrill}
\address[Ian Kiming]{Dept. of Mathematics, University of Copenhagen, Universitetsparken 5, 2100 Copenhagen \O , Denmark}
\email{kiming@math.ku.dk}
\address[Helena A. Verrill]{Dept. of Mathematics, Louisiana State University, Baton Rouge, Louisiana 70803-4918, USA}
\email{verrill@math.lsu.edu}

\begin{abstract} We give a parametrization of the possible Serre invariants $(N,k,\nu)$ of modular mod $\ell$ Galois
representations of the exceptional types $A_4$, $S_4$, $A_5$, in terms of local data attached to the fields cut out by the
associated projective representations. We show how this result combined with certain global considerations leads to an effective
procedure that will determine for a given eigenform $f$ and prime $\ell$ whether a mod $\ell$ representation attached to $f$ is
exceptional. We illustrate with numerical examples.
\end{abstract}
\maketitle

\section{Introduction.} \label{intro} Suppose that $f$ is a eigenform of weight $\ge$ 2 for some $\Gamma_1(M)$, and let $\ell$ be a
prime number. If $\lambda$ is a prime of $\overline\Q$ above $\ell$ then by a construction \cite{del} of Deligne followed by
reduction mod $\lambda$ and semisimplification, there is attached to $f$ a mod $\lambda$ representation $\rho \colon G_{\Q}
\longrightarrow \GL_2(\overline\F_{\ell})$ of the absolute Galois group of $\Q$. If $\rho$ is irreducible, the classification
\cite{dic} of finite subgroups of $\GL_2(\overline\F_{\ell})$ implies that the image of the projectivisation
$$
\Proj(\rho) \colon G_{\Q} \longrightarrow \PGL_2(\overline\F_{\ell})
$$
of $\rho$ is {\it (i)} dihedral, {\it (ii)} isomorphic to one of the groups $A_4$, $S_4$, $A_5$, or {\it (iii)} isomorphic to
$\PGL_2(\F_{\ell^n})$ or $\PSL_2(\F_{\ell^n})$ for some $n$. Thanks to work by Swinnerton--Dyer, Serre, Ribet and Momose, cf.
\cite{sw1}, \cite{sw2}, \cite{ser3}, \cite{rib2}, \cite{rib3}, \cite{mom}, one knows that the cases {\it (i)} and {\it (ii)}
occur for only a finite number of primes $\ell$. Also, if we are in case {\it (i)} or {\it (ii)} for our particular prime
$\lambda$, then this implies the existence of certain exceptional mod $\lambda$ congruences for the Fourier coefficients of $f$.
The link to exceptional mod $\lambda$ congruences is one the the main reasons for interest in determining for the given pair
$(f,\lambda)$ whether $\rho$ has `small image', i.e., whether we are in one of the cases {\it (i)} or {\it (ii)}.
\smallskip

It is thus natural to ask whether we can devise a method that will determine by a finite amount of computation the image of
$\rho$, or at least that of $\Proj(\rho)$. In principle the explicit version of Chebotarev's theorem \cite{lom} provides us with
such a method, but it seems impractical to use this since the bounds involved are generally astronomical. Of course, if one
expects to be in case {\it (iii)}, one can hope to quickly notice this by finding \dash--- via the Fourier coefficients of $f$
\dash--- elements in the image of $\rho$ whose presence imply that case {\it (iii)} prevails. This seems to work remarkably well,
see for example \cite{dv}. On the other hand, since case {\it (i)} means that $\rho$ is induced from a 1-dimensional character on
$G_M$ for some quadratic $M/\Q$, it is not hard to see how one can obtain through algebraic number theory combined with
\cite{stu} a method that will confirm or deny that case {\it (i)} prevails. These considerations may convince us that the cases
{\it (ii)} \dash--- which we refer to as the exceptional cases \dash--- are especially problematic to confirm for a given
$(f,\lambda)$. For example, in the works \cite{sw1}, \cite{sw2}, \cite{ser3} that deal with classical forms $f$ of level 1, it
was precisely the case where $f=Q\Delta$ is the normalized cusp form of weight 16 and level 1, $\ell=59$ and where $\Ima
\Proj(\rho)$ turns out to be of $S_4$-type, that caused the most trouble. This example was finally settled by Haberland in
\cite{hab} who worked hard with the Fourier coefficients of $Q\Delta$.
\smallskip

The contribution of the present article is the systematic development of the following simple idea: If $\rho$ is of the
exceptional type {\it (ii)}, then $\Proj(\rho)$ cuts out an extension $K/\Q$ with Galois group isomorphic to one of the groups
$A_4$, $S_4$, $A_5$. Then $K/\Q$ gives rise to a {\it complex} projective Galois representation $\pi$ (by choosing an embedding
$\Gal (K/\Q) \hookrightarrow \PGL_2(\C)$) which \dash--- in an appropriate sense \dash--- may be assumed to yield $\Proj(\rho)$
upon reducing mod $\lambda$, cf. the discussion in section \ref{lift} below. Then, if we lift $\pi$ to a linear representation $r
\colon G_{\Q} \rightarrow \GL_2(\C)$ and make this $\ell$-adically integral, $\rho$ is a twist of the mod $\lambda$ reduction of
$r$. Utilizing then the extensive information \dash--- as given in \cite{buh}, \cite{zin}, \cite{kim} \dash--- about the behavior
of complex lifts upon twisting, we are able to determine purely in terms of local data attached to the extension $K/\Q$ the
(infinitely many) possible Serre invariants $(N,k,\nu)$ of representations $\rho$ that cut out $K/\Q$ projectively. Here, and in
the following, by `Serre invariant' $(N,k,\nu)$ we mean the quantities attached by Serre in \cite{ser1} to an irreducible mod
$\ell$ representation $\rho$, with the small modification in connection with the weight $k$ as given in \cite{edi}. Combining
this with results on Serre's refined conjecture as given in \cite{rib}, \cite{dia}, \cite{edi}, \cite{cv}, we can predict the
minimal types $(N,k,\nu)$ of eigenforms whose attached projective mod $\ell$ representations cut out a given $K/\Q$. This is what
Theorem \ref{t1} below does: It gives a parametrization of these (infinitely many) minimal types $(N,k,\nu)$ purely in terms of
local data attached to $K/\Q$.

One remark about the theorem: We discuss only minimal types $(N,k,\nu)$ where every prime divisor of $N$ is a ramification point
for the given $K/\Q$; this is because, as one immediately checks, any other minimal type will then have the shape $(N\cdot
m^2,k,\nu \cdot \phi^2)$ where $\phi$ is a character of conductor $m$ whose prime divisors are all different from $\ell$ and
unramified in $K$.
\medskip

The information obtained from Theorem \ref{t1} is sufficient to actually prove exceptionality for some particular examples of
modular mod $\ell$ representations. We give some of these immediately in section \ref{exam} below. For instance, we are able to
quickly reprove Haberland's result in connection with $Q\Delta$ and $\ell =59$ without knowing a single Fourier coefficient of
the form. Another and somewhat more complicated example concerns the unique normalized cusp form of weight 4 on $\Gamma_0(8)$ and
$\ell =11$. These examples have also been dealt with by Boylan \cite{boy} by different methods.

We also discuss an $A_5$-type example occurring in Ribet's paper \cite{rib1}: Here we are able to immediately confirm \dash---
under an suitable modularity condition \dash--- a conjecture in that paper.
\medskip

In section \ref{lift} we discuss the above lifting of exceptional projective mod $\lambda$ representations to complex ones. After
that, in the beginning of section \ref{main} we describe first the structure of the proof of Theorem \ref{t1} in detail, and then
proceed with the actual proof. Finally, in section \ref{rem} we show how Theorem \ref{t1} can be complemented by certain global
considerations to yield an effective general procedure that will decide whether the exceptional case prevails for any explicitly
given pair $(f,\lambda)$.
\medskip

\noindent {\bf Notation.} The Galois group of a separable closure of a field $k$ will be denoted $G_k$.

The symbol $\ell$ denotes a fixed {\it odd} prime number, but $p$ may be any prime. We fix a prime $\lambda$ of $\overline\Q$
above $\ell$. More generally, if $p$ is any prime it will be convenient to assume that a prime of $\overline{\Q}$ over $p$ has
been fixed. Corresponding to this prime over $p$ we have a decomposition subgroup and an inertia subgroup of $G_{\Q}$; the
inertia subgroup will be denoted by $I_p$.

Generally, $K/{\Q}$ will denote a not totally real extension with Galois group isomorphic to one of the groups $A_4$, $S_4$,
$A_5$; in this situation, $\mathfrak l$ and $\p$ will denote primes of $K$ over $\ell$ and $p$ respectively.

All Galois representations occurring are assumed to cut out a finite extension of the base field in question.

If $\rho$ is a linear or projective representation of $G_{\Q}$, we shall denote by $\rho_p$ the restriction of $\rho$ to the
decomposition group in $G_{\Q}$ corresponding to the fixed prime of $\overline\Q$ over $p$. We then view $\rho_p$ as a
representation of $G_{\Q_p}$. By the {\it projective kernel field} of $\rho$ we mean the field fixed by the kernel of $\rho$ if
$\rho$ is projective, or the kernel of the projectivisation of $\rho$ if $\rho$ is linear.

The projectivisation of a given linear representation $\rho$ will be denoted by $\Proj(\rho)$.

We denote by $U_p$ the group of units in $\Z_p$. If $\e_p$ is a character of $\Q_p$, we denote by $c_p(\e_p)$ the {\it exponent}
of the conductor of $\e_p$.

We use (local) class field theory freely and without special notice. Thus for example, if $\e_p$ is a character of $G_{\Q_p}$ we
may view the restriction of $\e_p$ to the inertia subgroup of $G_{\Q_p}$ as a character on $U_p$.

The symbol $\chi$ always denotes the mod $\ell$ cyclotomic character $\chi \colon G_{\Q} \longrightarrow \F_{\ell}^{\times}$.

\section{Statement of results.}\label{results} Let $\ell$ be an odd prime, and suppose that $K/{\Q}$ is a non-real Galois
extension with Galois group $G$ isomorphic to one of the groups $A_4$, $S_4$, $A_5$.
\smallskip

Suppose that $p$ is a prime that ramifies in $K$ and that $\phi_p \colon U_p\rightarrow \C^{\times}$ is a homomorphism.
\smallskip

In section \ref{local} we will define a character
$$
\e _p = \e (K_{\p}/{\Q}_p) \colon ~U_p \longrightarrow \overline{\F}_{\ell}^{\times}
$$
depending only on the extension $K_{\p}/{\Q}_p$, and also non-negative integers
$$
n_p(\phi_p) = n_p(K_{\p}/{\Q}_p,\phi_p) \quad \mbox{and} \quad \delta_p(\phi_p) = \delta_p(K_{\p}/{\Q}_p,\phi_p)
$$
that depend on $K_{\p}/{\Q}_p$ as well as on $\phi_p$.
\smallskip

Also, in sections \ref{weight}--\ref{end} we define a set
$$
\W_{\ell} = \W (K_{\mathfrak l}/{\Q}_{\ell})
$$
of natural numbers depending on the extension $K_{\mathfrak l}/{\Q}_{\ell}$. The set $\W_{\ell}$ also comes equipped with a
partition $\W_{\ell} = \W_{\ell}^{+} \cup \W_{\ell}^{-}$ in case $\ell \not= 5$ and $K_{\mathfrak l}/{\Q}_{\ell}$ has degree
divisible by 5.
\smallskip

In this situation and with these definitions our main result can be formulated as follows.

\begin{thm} \label{t1} Retaining the above notation, let $S$ denote the product of the finite primes different from $\ell$
that ramify in $K$.

Suppose that $N$ is a natural number whose prime divisors all divide $S$, that
$$
\nu \colon (\Z/\Z N)^{\times} \longrightarrow \overline{\F}_{\ell}^{\times}
$$
is a homomorphism, and that $k\in \{ 2,\ldots ,\ell-1 \}$.
\medskip

(a) Suppose that $(N,k,\nu)$ is the triple of Serre invariants of a modular mod $\lambda$ representation cutting out $K/\Q$
projectively.
\smallskip

Then $k\in \W _{\ell}$ and there exists a global character of order prime to $\ell$
$$
\phi \colon G_{\Q} \rightarrow \C^{\times}
$$
unramified outside $S\cdot \infty$, such that
$$
N =\prod_{p\mid S} p^{n_p(\phi_{\mid I_p}) - \delta_p(\phi_{\mid I_p})} ~,
$$
and such that
$$
\nu(r) = \prod_{p\mid S} \e _p(r)^{-1} \cdot \left( \phi_p(r)^{-2} \mod \lambda \right) , \quad \mbox{for } (r,N)=1.
$$

(b) Conversely, suppose that $N$ and $\nu$ have the forms as in (a) for some global character $\phi$ of order prime to $\ell$.
Assume also that $K/\Q$ is cut out projectively by some ($2$-dimensional) modular mod $\lambda$ representation.
\smallskip

Then if either $K_{\mathfrak l}/{\Q}_{\ell}$ has degree not divisible by 5 or if $\ell=5$, there exists for any $k\in \W_{\ell}$
an eigenform of type $(N,k,\nu)$ and with $K$ the projective kernel field of the attached mod $\lambda$ representation.
\smallskip

If on the other hand, $K_{\mathfrak l}/{\Q}_{\ell}$ does have degree divisible by 5 and $\ell \not= 5$, there exists a $\mu \in
\{ +,-\}$ such that: For any $k\in \W_{\ell}^{\mu}$ there exists an eigenform of type $(N,k,\nu)$ and with $K$ the projective
kernel field of the attached mod $\lambda$ representation.
\end{thm}
\bigskip

For an outline of the proof of the theorem see the beginning of section \ref{main} below.
\bigskip

\noindent {\bf Remarks:} (1) In the theorem, the numbers $n_p(\phi_{\mid I_p})$ actually depend only on ($K_{\p}/\Q_p$ and) the
conductor of $\phi_p$ except possibly when $K_{\p}/\Q_p$ is cyclic.
\smallskip

\noindent (2) The small indeterminacy concerning the weight $k$ in part (b) of the theorem stems from the from the fact that
there are $2$ inequivalent embeddings $A_5\hookrightarrow \PGL_2(\overline{\F}_{\ell})$ if $\ell \not= 5$. It is true \dash--- in
an appropriate sense \dash--- that swapping between these 2 embeddings will change the sign $\mu$. The point is however that our
proof of the theorem is by purely local means and that the {\it simultaneous} specification of all of the 3 quantities $N$, $k$
and $\nu$ really requires some sort of global information in this case. An example of this type of phenomenon can be found in
\cite{buh}, Chap. 6. We shall indicate in section \ref{rem} below how the question can be resolved by certain global information
in any explicitly given case.
\smallskip

\noindent (3) By results of Langlands and Tunnell, cf. \cite{lan}, \cite{tun}, the modularity assumption in part (b) of the
theorem, i.e., that $K/\Q$ be cut out projectively by some modular mod $\lambda$ representation, is always true if $K/\Q$ is of
type $S_4$ or $A_4$. In case that $K/\Q$ is of $A_5$-type, the assumption is not an extremely strong one, as follows from recent
progress on the Artin conjecture for odd icosahedral representations of $G_{\Q}$, cf. \cite{bdst}, \cite{tay}.

One could define a non-real extension $K/\Q$ of one of the types $A_4$, $S_4$, or $A_5$ to be modular if it is the
projectivisation of the reduction mod $\lambda$ of the $2$-dimensional $\lambda$-adic representation attached by Deligne and
Serre \cite{ds} to an eigenform of weight $1$ for some $\Gamma_1(N)$. With this definition it would not be hard to show that
$K/\Q$ is modular if and only if $K/\Q$ is cut out by {\it some} modular mod $\lambda$ representation if and only if {\it any}
mod $\lambda$ representation cutting out $K/\Q$ projectively is modular (for this one needs the results of \cite{edi}, and in
case that $K/\Q$ is of $A_5$-type and $\ell \not= 5$ where we have $2$ inequivalent embeddings $A_5\hookrightarrow
\PGL_2(\overline{\F}_{\ell})$ as mentioned above, the argument requires a small special consideration involving the complex
conjugation of modular forms).
\smallskip

\noindent (4) It would have been possible to include a discussion of non-optimal levels (\cite{dt} and \cite{buz}). Also, we
could have discussed weights outside the interval $2\le k\le \ell-1$; the main purpose for restricting to this interval is to
exclude the case where the weight of the mod $\ell$ representation is 1. This case would occur if the given extension $K/\Q$ was
unramified over $\ell$ which is a rather uninteresting situation in the context of the present paper. Of course, for a general
modular mod $\ell$ representation the natural interval of weights to discuss is the full range $1\le k\le \ell+1$; however, for
representations of the exceptional types studied in this paper, restricting to the interval $2\le k\le \ell-1$ causes only a
small inconvenience in the very special case that $\ell=3$ and $K_{\mathfrak l}/\Q_{\ell}$ is totally ramified dihedral of order
6, in which case the weight $k=4$ may occur `naturally'. See the proof of part (b) of the Theorem for more details.

Finally, one should notice that the proof of the theorem obviously enables one to actually count the number of modular mod $\ell$
representations with a given triple of Serre invariants $(N,k,\nu)$ (of course under the appropriate modularity assumption in the
$A_5$-case and with the minor indeterminacy of part (b)).

These possible extensions of Theorem \ref{t1} have been discarded mainly in order to prevent overloading in the statement of the
theorem.

\section{Examples.} \label{exam} In this section we give a few simple examples of the use of Theorem \ref{t1} in connection with
questions about the image of modular mod $\ell$ representations.

\subsection{} \label{S4} An interesting class of $S_4$-type representations arises in connecting with non-real $S_4$-extensions
$K/\Q$ where $K_{\mathfrak l}/\Q_{\ell}$ is dihedral of order 8. In that case, we must have $\ell \equiv 3 \pod{4}$ as there are
no dihedral extensions of $\Q_{\ell}$ of order 8 if $\ell \equiv 1 \pod{4}$. Certainly, the extension $K/\Q$ is cut out
projectively by some modular mod $\lambda$ representation: Embedding $\Gal (K/\Q) \hookrightarrow \PGL_2(\C)$ and lifting this to
a linear representation, which is possible by a theorem of Tate, cf. \cite{ser2}, results of Langlands and Tunnell, \cite{lan},
\cite{tun}, imply that $K/\Q$ is cut out projectively by the complex representation attached to a modular form of weight 1.
Making this representation $\lambda$-adically integral and reducing mod $\lambda$ proves the claim.
\smallskip

Now, the definition of $\W_{\ell}$ in this case is:
$$
\W_{\ell} := \left\{ \frac{\ell+5}{4} , \frac{3\ell+7}{4} \right\} ~,
$$
cf. section \ref{tam_ram_dih} below. Part $(b)$ of Theorem \ref{t1} then informs us that there are exceptional modular mod
$\lambda$ representations with $K/\Q$ the projective kernel field for any of the weights $k\in \W_{\ell}$.
\smallskip

Let us now specialize even further to the case where $K/\Q$ is unramified outside $2\cdot \ell \cdot \infty$, and is either
unramified or of $S_4$-type locally above $2$. So, if $\p_2$ a prime of $K$ over $2$ we have either $K_{\p_2} = \Q_2$, or else
\dash--- see Weil \cite{wei} \dash--- $K_{\p_2}$ is one the the following three extensions $M_i/\Q_2$, $i=2,3,4$:
$$
M_i := \Q_2( \zeta , \omega , \sqrt{x_i} , \sqrt{\sigma x_i} , \sqrt{\sigma^2 x_i} )
$$
where $\zeta$ is a primitive 3'rd root of unity, $\omega^3 = 2$, $\sigma$ is the automorphism of $\Q_2(\zeta,\omega)$ given by
$\sigma \omega = \zeta \omega$, and $x_2 := (1+\omega)(1+\omega^2)(1+\omega^3)$, $x_3:= (1+\omega)(1+\omega^3)$, $x_4 :=
(1+\omega^2)$.
\smallskip

We have then the following immediate corollary to Theorem \ref{t1}.

\begin{cor}\label{c1} Suppose that $\ell \equiv 3 \pod{4}$ and that $K/\Q$ is a non-real $S_4$-type extension which is unramified
outside $2\cdot \ell \cdot \infty$, and such that $K_{\mathfrak l}/\Q_{\ell}$ is the unique dihedral extension of order $8$ of
$\Q_{\ell}$. Suppose further that $K$ is either unramified or of $S_4$-type locally above $2$. Let $\p_2$ a prime of $K$ over
$2$.

Let $\phi \colon \Z/\Z 2^c \rightarrow \C^{\times}$ be a Dirichlet character of conductor $2^c$, and define the non-negative
integer $n$ and the character $\e \colon ~U_2 \longrightarrow \overline{\F}_{\ell}^{\times}$ as follows:
$$
\left\{ \begin{array}{llll} \e = 1 ~~ & \mbox{and} ~~ & n = 2c & \mbox{if} ~~ K_{\p_2} = \Q_2 ~, \\
&&& \\
\e(-1) = -1,~ \e(5) = 1 ~~ & \mbox{and} ~~ & \left\{ \begin{array}{ll} n = 7 & \mbox{if} ~~ c\le 3 \\ n = 2c & \mbox{if} ~~
c\ge 4 \end{array} \right\} & \mbox{if} ~~ K_{\p_2} = M_2 ~, \\
&&& \\
\e = 1 ~~ & \mbox{and} ~~ & \left\{ \begin{array}{ll} n = 7 & \mbox{if} ~~ c\le 3 \\ n = 2c & \mbox{if} ~~
c\ge 4 \end{array} \right\} & \mbox{if} ~~ K_{\p_2} = M_3 ~, \\
&&& \\
\e = 1 ~~ & \mbox{and} ~~ & \left\{ \begin{array}{ll} n = 3 & \mbox{if} ~~ c\le 1 \\ n = 2c & \mbox{if} ~~ c\ge 2 \end{array}
\right\} & \mbox{if} ~~ K_{\p_2} = M_4 ~.
\end{array} \right.
$$

Then for each of the weights $k=\frac{\ell+5}{4} , \frac{3\ell+7}{4}$, there exists a mod $\lambda$ eigenform whose attached mod
$\lambda$ representation cuts out the extension $K/\Q$ projectively, and which has type $(2^n,k,\nu)$ with $\nu$ given by:
$$
\nu(r) =  \e(r)^{-1} \cdot \left( \phi(r)^{-2} \mod \lambda \right) , \quad \mbox{for odd } r.
$$
\end{cor}
\begin{proof} The result follows immediately from Theorem \ref{t1} once one reviews for this particular case the definitions, in
section \ref{primitive}, \ref{deg}, and \ref{tam_ram_dih} below, of the numbers $n_2(\phi_{\mid I_2})$ and $\delta_2(\phi_{\mid
I_2})$, of the character $\e = \e_2 \colon ~ U_2 \rightarrow \overline{\F}_{\ell}^{\times}$, and of the set $\W_{\ell}$.
\end{proof}

\subsubsection{} Suppose additionally that $K/\Q$ is unramified outside $\ell \cdot \infty$. In \cite{dou} such fields were
studied using other methods than ours and under an additional condition on $\ell$ made to ensure the solvability of a certain
embedding problem associated with lifting $G_{\Q} \twoheadrightarrow \Gal (K/\Q) \hookrightarrow \PGL_2(\C)$ to a linear
representation. By solving explicitly this embedding problem and studying the behavior of the solution locally at $\ell$ the
author of that paper was able to find the connection to cusp forms of weight $\frac{\ell+5}{4}$. Here we proceed without any
further assumptions than those already made on $K$ and conclude from Theorem \ref{t1} that there are modular mod $\ell$
representations of level 1 and the above weights with $K$ the associated projective kernel field. The classical example of this
phenomenon concerns the unique cusp form $Q\Delta$ of level 1 and weight 16 in connection with $\ell=59$. This case was first
considered by Swinnerton--Dyer in \cite{sw1} who found strong reasons to believe that the mod 59 representation attached to
$Q\Delta$ is exceptional of type $S_4$. This was then subsequently proved in \cite{hab} by Haberland who \dash--- having less
powerful theorems at his disposal than we do \dash--- worked hard numerically with the Fourier expansion of $Q\Delta$. We can
reprove the statement easily and without knowing any Fourier coefficients: One knows (cf. \cite{hab}) that there is a unique
$S_4$-extension of $\Q$ unramified outside 59 and $\infty$. In fact, the splitting field $K$ of the polynomial
$$
x^4-x^3-7x^2+11x+3
$$
is seen to be such a field. It is non-real, and dihedral of order 8 over 59. Theorem \ref{t1} implies the existence of modular
mod 59 representations of level 1 and weights $\frac{59+5}{4}=16$, $\frac{3\cdot 59+7}{4}=46$ with $K$ as projective kernel
field. Since $Q\Delta$ is the unique cusp form of level 1 and weight 16, the claim follows.

\subsubsection{} A somewhat more complicated example occurs in connection with the mod 11 representation $\rho$ attached to
the unique cusp form $\eta(2z)^4\cdot \eta(4z)^4$ of weight 4 on $\Gamma_0(8)$ where $\eta(z)$ is the Dedekind eta function:
$$
\eta(z)=e^{2\pi iz/24}\prod_{n=1}^{\infty}(1-e^{2\pi inz}) ~.
$$
Numerical experiments performed by the second author suggested strongly that $\rho$ is exceptional of $S_4$-type; verifying this
was one of the initial starting points of the present article. If we wish to verify this via Theorem \ref{t1} we must begin by
looking for a candidate $K/\Q$ for the corresponding projective kernel field. Thus, $K/\Q$ should be a non-real $S_4$-extension
unramified outside $2\cdot 11\cdot \infty$. By class field theory one finds that the unique $S_3$-subextention $L/\Q$ contained
in $K/\Q$ would have to be the splitting field of
$$
x^3-x^2-x-1 ~,
$$
which is in fact an $S_3$-extension of $\Q$ unramified outside $2\cdot 11\cdot \infty$. Now we have to look for an embedding of
this $L/\Q$ in an $S_4$-extension $K/\Q$ of the desired type. It turns out that there is more than one such $K/\Q$ but we shall
focus on the one that is relevant for our example. We seek to find the candidate $K/\Q$ by utilizing the method of \cite{bf}:

Consider the elliptic curve $-22y^2=x^3-x^2-x-1$ which has minimal Weierstrass model
$$
E\colon ~ y^2 = x^3+x^2-645x+14771 ~.
$$
The curve $E$ has good reduction outside $2\cdot 11$, and has the rational point
$$
P=(-26,121) \in E(\Q)\backslash 2\cdot E(\Q) ~.
$$
According to \cite{bf}, the field generated over $\Q$ by the $x$-coordinates of points $Q\in E(\overline\Q)$ with $2\cdot Q=P$ is
a (non-real) $S_4$-extension $K/\Q$ unramified outside $2\cdot 11\cdot \infty$. Using the 2-division polynomial for $E$ we find
that this $K$ is the splitting field of the polynomial
$$
x^4+104x^3+1394x^2-185248x+1893125 ~.
$$
Applying to this polynomial the function {\tt polred} in PARI, cf. \cite{pari}, we find that $K$ can also be given as the
splitting field of
$$
x^4-2x^3-4x^2+16x-13 ~.
$$
One can then verify directly that $K/\Q$ is indeed a non-real $S_4$-extension unramified outside $2\cdot 11\cdot \infty$, that
the local extension over 11 is the unique $D_4$-extension of $\Q_{11}$, and finally that the local extension over 2 is the
$S_4$-extension that was denoted $M_4$ in \cite{kim} and in \ref{S4} above. Thus, Corollary \ref{c1} applies and gives the
existence of an exceptional $S_4$-type modular mod 11 representation of level 8, weight 4 and trivial nebentypus, which must then
be our representation $\rho$.
\smallskip

Similarly, one can show that the mod 19 representation attached to the unique cusp form $\eta(2z)^{12}$ of weight 6 on
$\Gamma_0(4)$ is exceptional of $S_4$-type with projective kernel field the splitting field of
$$
x^4-x^3-2x^2-6x-2 ~.
$$

\subsection{} \label{A5} An interesting $A_5$-type example occurs in Ribet's paper \cite{rib1}, p. 284: There are exactly 2
newforms of weight 2 on $\Gamma_0(23)$; call them $f_1,f_2$. The example concerns the mod 3 representations $\rho_i$ attached to
$f_i$. Arguments given {\it loc. cit.}, suggested that each of the representations $\rho_i$ is of $A_5$-type, and in fact with
projective kernel field the $A_5$-extension of $\Q$ given as the splitting field $K$ of the polynomial
$$
x^5+3x^3+6x^2+9 ~.\leqno{(\ast)}
$$
It was further stated {\it loc. cit.}, that computations performed by Mestre were consistent with this claim. We shall show here
how Theorem \ref{t1} immediately implies this statement under the assumption of modularity of a 2-dimensional complex
representation associated with $K$.

The polynomial $(\ast)$ occurs as the second entry of the table in \cite{bk}, and if we choose an embedding $A_5\hookrightarrow
\PGL_2(\C)$ we obtain a complex, odd $A_5$-type representation $\rho$ of conductor $3^3\cdot 23^2$. It should be fairly clear
that a rigorous verification of the above statement about the $\rho_i$ must involve either a proof of the Artin conjecture of
this complex $\rho$ or else the application of some form of explicit Chebotarev. The authors have not considered how to
practically implement the second option. Concerning the first, let us remark that the results of the papers \cite{bdst} and
\cite{tay} do not suffice to prove the Artin conjecture for $\rho$ (in $K$, a Frobenius over $2$ has order $5$, and $3$ is
totally ramified of degree $6$). Furthermore, an attempt at a direct computational verification of Artin's conjecture for this
$\rho$ utilizing methods similar to those employed in \cite{kw}, can be estimated to lead to an intractably large computation.
Thus, we shall be content with giving a relative statement: Assuming the weight $1$ modularity for this particular $\rho$, the
$\rho_i$ are both of $A_5$-type with projective kernel field the above $K$. This is readily done: Denoting by $\p_{23}$ and
$\p_3$ primes of $K$ over 23 and 3 respectively, one finds that $K_{\p _{23}}/{\Q}_{23}$ is the unique dihedral extension of
degree 6 of ${\Q}_{23}$, and that $K_{\p_3}/{\Q}_3$ is totally ramified dihedral of degree 6. One checks that the latter
extension is peu ramifi\'{e}: An easy calculation shows that peu and tr\`{e}s ramifi\'{e} implies the contribution $3^4$ and
$3^6$ respectively to the discriminant of a quintic subfield of $K$ (cf. also Table 3.1 of \cite{buh}); as this discriminant is
in fact $3^4\cdot 23^2$ (\cite{bk}, Table 1), $K_{\p_3}/{\Q}_3$ is peu ramifi\'{e}.

The recipes in section \ref{main} then give: $n_{23}(1) = 1$, $\e_{23}=1$, and $\W_3 = \{ 2\}$. Under the above weight $1$
modularity assumption Theorem \ref{t1} thus implies the existence of a modular mod $3$ representation $\rho$ of level 23, weight
2 and trivial character that cuts out this $K/\Q$ projectively. On the other hand, choosing the other embedding
$A_5\hookrightarrow \PGL_2(\C)$, we obtain similarly another modular mod 3 representation which is not equivalent to $\rho$ but
has the same level, weight and character. As $S_2(\Gamma_0(23))$ has dimension $2$, the desired conclusions about the $\rho_i$
follow.

\section{Lifting exceptional mod \texorpdfstring{$\ell$}{l} representations to complex representations.}
\label{lift} In this section we prove some simple statements about lifting mod $\ell$ representations to complex ones. Everything
in this section is probably fairly standard knowledge but we have included it for lack of appropriate references.
\smallskip

Though we need only the case $n=2$, we start by discussing the more general case of $n$-dimensional representations as this
causes no extra work.
\smallskip

Given an irreducible Galois representation $\pi \colon G_{\Q} \rightarrow \PGL _n({\C})$ with finite image, by the reduction mod
$\lambda$ of $\pi$ we mean the representation $\Proj(\bar r) \colon G_{\Q} \rightarrow \PGL_n(\overline{\F}_\ell)$ constructed as
follows:

Let $r \colon G_{\Q} \rightarrow \GL _n({\C})$ be a lift of $\pi$, i.e., a representation which has projectivisation equal to
$\pi$. By a theorem of Tate (cf. \cite{ser2}) such a lift $r$ does exist. We may consider $r$ as a representation $r \colon
G_{\Q}\rightarrow \GL_n(\overline{\Q})$. Denote by $\overline{\Z}$ the ring of integers of $\overline{\Q}$. Since every finitely
generated ideal of $\overline{\Z}$ is principal, we can conjugate $r$ to a representation $G_{\Q} \rightarrow
\GL_n(\overline{\Z})$ which by abuse of notation we still denote by $r$. Now let $\bar r$ be the mod $\lambda$ reduction of $r$,
and consider the projectivisation $\Proj(\bar r) \colon G_{\Q} \longrightarrow \PGL_n(\overline\F_{\ell})$ of $\bar r$. One
easily checks that $\Proj(\bar r)$ is independent up to equivalence of the various choices made. We refer then to $\Proj(\bar r)$
as the mod $\lambda$ reduction of $\pi$. For clarity, the four representations $\pi$, $r$, $\bar r$, and $\Proj(\bar r)$ are
shown in the following diagram:
$$
\xymatrix{ &G_{\Q}\ar[dl]_r \ar[dr]^{\pi} \ar[ddl]^{\bar r} \ar[ddr]_{\Proj(\bar r)} \\
\GL_n(\overline{\Z})\ar[d] && \PGL_n(\overline{\Z})\ar[d] \\
\GL_n(\overline{\F}_{\ell})\ar[rr] && \PGL_n(\overline{\F}_{\ell}) }
$$

Reversing the roles of $\pi$ and $\Proj(\bar r)$, if $\Pi \colon G_{\Q} \rightarrow \PGL _n(\overline{\F}_\ell)$ is given and
assumed to be irreducible, we may ask for a representation $\pi \colon G_{\Q} \rightarrow \PGL _n({\C})$ whose mod $\lambda$
reduction is the given $\Pi$. If such a $\pi$ exists we may and will then refer to it as a {\it complex lift} of $\Pi$.

\begin{prop} \label{p1} Let $\Pi \colon G_{\Q} \rightarrow \PGL _2(\overline{\F}_\ell)$ be an irreducible
representation such that $\Ima \Pi$ embeds in $\PGL_2(\C)$. Then $\Pi$ has a complex lift.
\end{prop}

We start the proof of the proposition by an observation that is best isolated as a small lemma:

\begin{lem} \label{l1} Let $\Pi \colon G_{\Q} \rightarrow \PGL _n(\overline{\F}_\ell)$ be an irreducible
representation such that the finite group $\Ima \Pi$ embeds in $\PGL _n(\C)$, and is $\ell$-solvable. Then $\Pi$ has a complex
lift.
\end{lem}

\begin{proof} Let us first notice that we have
$$
\Hc^2(G_{\Q},\overline{\F}_\ell^\times ) = 0 ~,
$$
where the action of $G_{\Q}$ on $\overline{\F}_\ell^\times$ is trivial. This is because we have
$$
\overline{\F}_\ell^\times \cong \prod_{p\not = \ell} {\Q}_p/{\Z}_p ~,
$$
and because we know that
$$
\Hc^2(G_{\Q},{\Q}_p/{\Z}_p) = 0 ~,
$$
cf. \cite{ser2}, \S 6.5. Consequently, $\Pi$ has some lift $\rho \colon G_{\Q} \rightarrow \GL_2(\overline{\F}_\ell)$. Now, since
$\Ima \Pi$ is $\ell$-solvable then so is the image of $\rho$. So we may apply the theorem of Fong--Swan, cf. for example
Theor\`{e}me 3, p. III-19 of \cite{ser4}, to infer that $\rho$ is in fact in the above sense the mod $\lambda$ reduction of a
complex representation $r$. The projectivisation of any such $r$ is then a complex lift of $\Pi$.
\end{proof}

\begin{proof}[Proof of Proposition \ref{p1}:] By the lemma we may assume that $G:=\Ima \Pi$ is not
$\ell$-solvable. As $G$ embeds into $\PGL_2(\C)$ this implies $G\cong A_5$ and $\ell \in \{ 3,5\}$. For these cases we finish the
proof by explicitly displaying a lift:
\smallskip

We first choose $m$ large enough so that $\Pi$ can be realized over $\F_{\ell^m}$, and so that $\F_{\ell^m}$ contains a primitive
5'th root of unity if $\ell \not =5$, i.e., if $\ell=3$. Now recall the well-known fact that there is at most one conjugacy class
of subgroups isomorphic to $G$ in $\PGL_2(\F_{\ell^m})$; this follows from the discussions in \cite{dic}, \S \S 258--259.

Consequently, any representation $G_{\Q} \longrightarrow \PGL_2(\F_{\ell^m})$ with the same kernel field as $\Pi$ will \dash---
up to equivalence \dash--- have the shape $\alpha \circ \Pi$ where $\alpha$ is some automorphism of $G=\Ima \Pi \le
\PGL_2(\F_{\ell^m})$. If $\ell=5$ then $\alpha$ must in fact be conjugation by some element of $\PGL_2(\F_{\ell^m})$, since
$G\cong A_5\cong \PSL_2({\F}_5)$ has automorphism group $\PGL_2({\F}_5)$ acting via conjugation (\cite{svdw}). On the other hand,
if $\ell=3$ there are precisely 2 inequivalent embeddings $G\hookrightarrow \PGL_2(\F_{\ell^m})$, as $\PGL_2(\F_{\ell^m})$ will
not in this case contain a subgroup isomorphic to $\PGL_2({\F}_5)$ (and as $\PSL_2({\F}_5)$ has index 2 in $\PGL_2({\F}_5)$).
\smallskip

In each of the cases we are considering, we now display an explicit complex lift of a given embedding $G\hookrightarrow
\PGL_2(\F_{\ell^m})$. In terms of the presentation $G\cong A_5 = \left\langle x,y \mid ~x^2=y^5=(xy)^3=1 \right\rangle$, the
following assignments define an embedding $G\hookrightarrow \PGL_2(\F_{\ell^m})$:
$$
x\mapsto \left( \begin{smallmatrix} -c & \omega \\
\omega & c \end{smallmatrix} \right),\quad y\mapsto
\left( \begin{smallmatrix} \varepsilon^2 & -\omega \\
0 & \varepsilon^{-2} \end{smallmatrix} \right),\leqno{(\sharp)}
$$
where $\omega:=\varepsilon+\varepsilon^{-1}$, $c:=\varepsilon^2-\varepsilon^{-2}$, and where $\varepsilon$ is a primitive 5'th
root of unity if $\ell=3$, but $\varepsilon =1$ if $\ell=5$. For $\ell=3$ the substitution $\varepsilon \mapsto \varepsilon^2$
brings one embedding to the other inequivalent one. Viewing the matrices in $(\sharp)$ as elements of $\GL_2({\C})$ by
interpreting $\varepsilon$ as a primitive complex 5'th root of unity, we have explicit complex lifts in all cases.
\end{proof}

\noindent {\bf Remark:} The authors do not know to what extent Proposition \ref{p1} generalizes to higher-dimensional situations.
The argument given for $n=2$ in case $\Ima \Pi$ is not $\ell$-solvable is admittedly awkward, relying as it does on the
classification of subgroups of $\PGL _2({\F}_{\ell^m})$. But we have not been able to find a more conceptual approach in this
case.

\section{Proof of Theorem \ref{t1}.} \label{main} In this section we prove Theorem \ref{t1}. So, let $K/{\Q}$, $G$ and $S$ be as
in the statement of the Theorem. Definitions of the quantities $n_p(K_{\p}/{\Q}_p,\phi_p)$ and $\e (K_{\p}/{\Q}_p)$ will be given
in section \ref{local} below, and the definition of the set $\W_{\ell} (K_{\mathfrak l}/{\Q}_{\ell})$ is stated in section
\ref{weight}.
\smallskip

We start with the proof of part (a). Let
$$
\rho \colon G_{\Q} \longrightarrow {\GL}_2(\overline{\F}_{\ell})
$$
be an irreducible modular mod $\lambda$ representation with Serre invariants $(N,k,\nu)$ where $2\le k\le \ell-1$, every prime
divisor of $N$ divides $S$, and such that the extension cut out by the projectivisation
$$
\Proj(\rho)\colon G_{\Q} \longrightarrow \PGL_2(\overline{\F}_{\ell})
$$
of $\rho$ is the given field $K$.
\smallskip

Let $\pi \colon G_{\Q} \longrightarrow \PGL_2(\C)$ denote a complex lift of $\Proj(\rho)$ in the sense of section \ref{lift}
(which exists according to Proposition \ref{p1}).
\smallskip

Now, if $r \colon G_{\Q} \longrightarrow {\GL}_2(\C)$ is any lift of $\pi$ to a linear representation unramified outside $S \ell
\cdot \infty$, and if $\bar r \colon G_{\Q} \longrightarrow {\GL}_2(\overline{\F}_{\ell})$ is its mod $\lambda$ reduction, the
definition of the notion `complex lift' implies that $\rho$ is a twist of $\bar r$ by some character $G_{\Q} \longrightarrow
\overline{\F}_{\ell}^{\times}$; we may write this character as $\bar{\phi} \cdot \chi^i$ for some $i\in \Z$ and with $\bar{\phi}$
a character unramified at $\ell$. In other words, if we make $r$ integral in the sense of the previous section, we have the
diagram:
$$
\xymatrix{ &G_{\Q}\ar[dl]_r \ar[dr]^{\pi} \ar[ddl]^{\bar r} \ar[ddr]_{\Proj(\rho)} \\
\GL_2(\overline{\Z})\ar[d] && \PGL_2(\overline{\Z})\ar[d] \\
\GL_2(\overline{\F}_{\ell})\ar[rr] && \PGL_2(\overline{\F}_{\ell}) }
$$
where each of the 3 `triangles' commute, and we have
$$
\rho = \bar r \otimes \bar{\phi} \chi^i ~.
$$

We may and will consider $\bar{\phi}$ as the mod $\lambda$ reduction of a complex character $\phi \colon G_{\Q} \longrightarrow
{\C}^{\times}$ of order prime to $\ell$, and view the restriction $\phi_p :=\phi_{\mid I_p}$ for $p\mid S$ as a homomorphism
$U_p\rightarrow \C^{\times}$.
\smallskip

In the following subsection we utilize the results of \cite{kim} to show the existence of a {\it particular} linear lift $r$ as
above which is unramified outside $S\ell \cdot \infty$ and has the property that each of the quantities
$$
c_p(r \otimes \phi) , \quad \mbox{for } p\mid S ~,
$$
and
$$
\det (r)_{\mid U_p}, \quad \mbox{for } p\mid S ~,
$$
can (and will) be explicitly determined. The first of these quantities depends on the data $K_{\p}/{\Q}_p$ and $\phi_p$, whereas
the second depends only on $K_{\p}/{\Q}_p$. Call the first of these quantities $n_p(\phi_p) = n_p(K_{\p}/{\Q}_p , \phi_p)$, and
let $\e_p = \e_p(K_{\p}/{\Q}_p)$ denote the mod $\lambda$ reduction of the second.
\smallskip

In section \ref{deg} we compute the number
$$
\delta_p(\phi_p) = \delta_p(K_{\p}/{\Q}_p,\phi_p) := n_p(\phi_p) - c_p(\bar r \otimes \bar{\phi})
$$
for $p\mid S$, which depends only on the datum $(K_{\p}/{\Q}_p,\phi_p)$. Then the conductor of $\rho = \bar r \otimes \bar\phi
\chi ^i$ is
$$
\prod_{p\mid S} p^{n_p(\phi_p) - \delta_p(\phi_p)} ~,
$$
which by the assumption on $\rho$ coincides with $N$. Further, for $p$ dividing $S$, the restriction of the determinant of $\rho$
to $I_p$ is the character $\prod_{p\mid S} \e_p \cdot \left( \phi_p^2 \mod \lambda \right)$ which by definition coincides with
$\nu_{\mid I_p}$ when $\nu$ is viewed, via global class field theory, as a character of $G_{\Q}$. As $\nu$ is unramified outside
$S\cdot \infty$ we have \dash--- again by global class field theory \dash--- that $\nu$ is given by
$$
x\mapsto \prod_{p\mid S} \e_p(x)^{-1} \cdot \left( \phi_p(x)^{-2} \mod \lambda \right) , \quad \mbox{for } (x,N)=1 ~,
$$
as a character on $(\Z/\Z N)^{\times}$.

In subsection \ref{weight} we study the local behavior of $\bar r \otimes \chi ^i$ at $\ell$ and define the set $\W _{\ell} = \W
(K_{\mathfrak l}/{\Q}_{\ell})$ such that $k\in \W _{\ell}$ may be inferred. This will then complete the proof of part (a) of
Theorem \ref{t1}.

\subsection{Local lifts: Conductors and determinants.}
\label{local} First we use a theorem of Tate, cf. \cite{ser2}, to reduce to a purely local question the problem of specifying a
convenient lift $r$ of $\pi$ to a linear representation: Define the representation $\pi_p \colon G_{\Q_p} \longrightarrow
\PGL_2(\C)$ as the restriction of $\pi$ to $G_{\Q_p} \le G_{\Q}$. Thus $\pi_p$ factors as
$$
\pi_p \colon G_{\Q_p} \twoheadrightarrow \Gal (K_{\p}/\Q_p) \hookrightarrow \PGL_2(\C) ~.
$$
Now, if we choose for each $p\mid S\cdot \ell$ a lift $r_p$ of $\pi_p$ to a linear representation, the above theorem of Tate
states that there is a uniquely determined lift $r$ of $\pi$ which is unramified outside $S\ell \cdot \infty$ and satisfies:
$$
(r)_{\mid I_p} \cong (r_p)_{\mid I_p} \quad \mbox{for } p\mid S\ell ~.
$$
We shall now address the problem of specifying for each $p$ a {\it particular} choice of $r_p$ with the property that the
quantities
$$
c_p(r_p \otimes \phi_p) \quad \mbox{and} \quad \det (r_p)_{\mid U_p} ~,
$$
can be explicitly given for any character $\phi_p \colon U_p\rightarrow \C^{\times}$. Then
$$
c_p(r \otimes \phi) = c_p(r_p \otimes \phi_{\mid I_p}) ~,
$$
and
$$
\det (r)_{\mid U_p} = \det (r_p)_{\mid U_p}
$$
are also explicitly known for any global character $\phi \colon G_{\Q} \rightarrow \C^{\times}$. This then completes the first
step of the proof of Theorem \ref{t1} as described above.
\medskip

The theory of local lifts $r_p$ was initiated by Weil in \cite{wei}, further developed by Buhler and Zink, cf. \cite{buh},
\cite{zin}, and complemented by the results in \cite{kim} where a full solution to the problem was given in all cases. Thus, we
can simply quote from the latter article and briefly review the results given there. There are 3 subcases corresponding to
whether the image of $\pi_p$ is cyclic, dihedral, or isomorphic to $A_4$ or $S_4$, respectively. Surprisingly, the most difficult
of these subcases is the dihedral one. For the following review, let $\phi_p$ denote an arbitrary character $G_{\Q_p} \rightarrow
\C^{\times}$.

\subsubsection{Cyclic cases.} This case is trivial: If $\Ima \pi_p$ is cyclic then $\pi_p$
is equivalent to
$$
\pi_p \sim \left( \begin{smallmatrix} \alpha &  0 \\  0 &  1 \end{smallmatrix} \right) \colon G_{\Q_p} \longrightarrow \PGL_2(\C)
~,
$$
for some character $\alpha \colon G_{\Q_p} \rightarrow \C^{\times}$ which we may and will view as a character on $\Q_p^{\times}$.
We can then choose $r_p \cong \alpha \oplus 1$ which gives
$$
n_p(\phi_p):=c_p(r_p \otimes \phi_p) = c_p(\alpha_{\mid U_p} \cdot \phi_p) + c_p(\phi_p) ~,
$$
and
$$
\e_p := (\alpha_{\mid U_p} \mod \lambda) ~.
$$

\subsubsection{Dihedral cases.} Suppose that $\Ima \pi_p$ is a dihedral group. Then $K_{\p}/\Q_p$
contains a quadratic extension $M/\Q_p$ such that $K_{\p}/M$ is cyclic and ramified. If the degree $[K_{\p}:\Q_p]$ is greater
than 4, the field $M$ is uniquely determined. If $[K_{\p}:\Q_p] = 4$ we let $M$ be any quadratic subfield of $K_{\p}$ such that
$K_{\p}/M$ is ramified. Let $\alpha \colon \Q_p^{\times} \rightarrow \C^{\times}$ be the quadratic character corresponding to $M$
via local class field theory. Now, the restriction of $\pi_p$ to $G_M$ has the form
$$
\left( \begin{smallmatrix} \beta & 0 \\ 0 & 1 \end{smallmatrix} \right)
$$
for some character $\beta$ on $G_M$ which by local class field theory may and will be viewed as a character on $M^{\times}$. In
Theorem 1 of \cite{kim} a special lift $r_p$ of $\pi_p$ is given, together with explicit formulas, depending on the data
$(M,\beta)$, for $\det (r_p)_{\mid U_p}$ and for $c_p(r_p \otimes \phi_p)$. So, we let
$$
\e_p := \left( \det (r_p)_{\mid U_p} \mod \lambda \right)
$$
and
$$
n_p(\phi_p) := c_p(r_p \otimes \phi_p)
$$
where $r_p$ is the special lift given loc. cit.; we have then explicit formulas for these quantities, but as these are somewhat
involved in the general case we shall not restate them in all detail but will be content to review their general shape and give a
few special cases:
\smallskip

Let $\omega$ be a uniformizer of $M$, let $b\ge 1$ be the exponent of the conductor of $\beta$ (i.e., $\beta$ has conductor
$\omega^b$), and let $t$ denote the break in the sequence of ramification groups of ${\mathfrak G} := \Gal (M/\Q_p)$ in the lower
numbering:
$$
{\mathfrak G} = {\mathfrak G}_0 = \ldots = {\mathfrak G}_t \not= {\mathfrak G}_{t+1} =0 ~.
$$
Then $c_p(r_p \otimes \phi_p)$ is given via:
$$
p^{n_p(\phi_p)} = p^{c_p(r_p \otimes \phi_p)} = D(M/\Q_p) \cdot N_{M/\Q_p}(\omega^{\max \{ t+b,\gamma(\phi_p) \} } ) ~,
$$
where $D(\cdot)$ and $N_{M/\Q_p}(\cdot)$ denote discriminant and norm respectively, and where
$$
\gamma(\phi_p)  ~\left\{ \begin{array}{ll} = c_p(\phi_p) & \mbox{ if } M/\Q_p \mbox{ unramified} \\
\le t+1 & \mbox{ if } M/\Q_p \mbox{ ramified and } c_p(\phi_p) \le t+1 \\ = 2\cdot c_p(\phi_p) +1-t & \mbox{ if } M/\Q_p \mbox{
ramified and } c_p(\phi_p) \ge t+2. \end{array} \right.
$$
Thus, as is immediately seen, the number $n_p(\phi_p)$ depends only on $K_{\p}/\Q_p$ and on the conductor of $\phi_p$ (rather
than on $\phi_p$ itself).
\smallskip

The character $\det (r_p)_{\mid U_p}$ has the general shape
$$
\det (r_p)_{\mid U_p} = \alpha_{\mid U_p} \cdot \psi_1 \cdot \psi_2
$$
where $\psi_1, \psi_2$ are the characters given loc. cit., which are both of $2$-power order. The character $\psi_1$ is at most
tamely ramified, whereas $\psi_2$ is wildly ramified if it is non-trivial, which may happen only if $p=2$. The definitions of
$\psi_1$ and $\psi_2$ are complicated in general and will not be restated, but we give a few special cases:
\smallskip

If the degree $m:=[K_{\p}:M]$ is odd then $\psi_1=\psi_2 =1$ so that $\e_p$ is the trivial or quadratic character
$$
\e_p = (\alpha_{\mid U_p} \mod \lambda) ~.
$$

If additionally we have $p\nmid m$ then one finds $n_p(1) =2$.
\smallskip

On the other hand, if $K_{\p}/\Q_p$ is dihedral of order 8, but $p\not = 2$, then $p\equiv 3 \pod{4}$ and $K_{\p}$ is unique. One
has in this case that $M$ is the unramified quadratic extension $\Q_p(\sqrt{-1})/\Q_p$ so that $\alpha$ is unramified. Further,
$\psi_2=1$, but $\psi_1$ is the non-trivial quadratic character on $\{ \pm 1\} \le \Q_p^{\times}$. With the above notation one
has $b=1$ and $t=0$ so that
$$
n_p(\phi_p) = 2\cdot \max \{ 1, c_p(\phi_p) \} ~.
$$

\subsubsection{The `primitive' cases for \texorpdfstring{$p=2$}{p=2}.}\label{primitive} Finally,
$\Ima \pi_p$ may be isomorphic to $A_4$ or $S_4$ in which case we must have $p=2$. The study of these cases was initiated by Weil
in \cite{wei}. In \cite{buh} Buhler computed the minimal possible conductor of a lift of $\pi_2$ and subsequently Zink determined
in \cite{zin} all possible conductors of lifts of $\pi_2$. In \cite{kim} these results were complemented by a discussion of the
associated determinants. We can refer to Theorem 2 of \cite{kim} where a special lift $r_2$ of $\pi_2$ was given together with
formulas for the restriction of its determinant to $U_2$ and for the conductor of $r_2 \otimes \phi_2$. We define then
$$
n_2(\phi_2) := c_2(r_2 \otimes \phi_2)
$$
and
$$
\e_2 := \left( \det (r_2)_{\mid U_2} \mod \lambda \right) ~,
$$
where $r_2$ is a lift of $\pi_2$ as in Theorem 2 of \cite{kim}. Thus the formulas of that theorem give \dash--- depending on the
extension $K_{\p_2}/\Q_2$ ($\p_2$ a prime of $K$ over 2) \dash--- explicit formulas for the quantities $n_2(\phi_2)$ and $\e_2$.
The character $\e_2$ is at most quadratic, and the number $n_2(\phi_2)$ actually only depends on $c_2(\phi_2)$ in every case.

According to \cite{wei} there is exactly $1$ extension of $\Q_2$ of $A_4$-type; call it $M_1$. Furthermore, there are exactly $3$
extensions of $S_4$-type; as already noted in section \ref{exam} above these are $M_i/\Q_2$, $i=2,3,4$, where:
$$
M_i := \Q_2( \zeta , \omega , \sqrt{x_i} , \sqrt{\sigma x_i} , \sqrt{\sigma^2 x_i} )
$$
where $\zeta$ is a primitive 3'rd root of unity, $\omega^3 = 2$, $\sigma$ is the automorphism of $\Q_2(\zeta,\omega)$ given by
$\sigma \omega = \zeta \omega$, and $x_2 := (1+\omega)(1+\omega^2)(1+\omega^3)$, $x_3:= (1+\omega)(1+\omega^3)$, $x_4 :=
(1+\omega^2)$.
\medskip

Theorem 2 of \cite{kim} then gives the above number $n_2(\phi_2)$ and the restriction of the above character $\e_2$ to $U_2$ as
follows: Write $n_2 := n_2(\phi_2)$ and $c_2 := c_2(\phi_2)$ for convenience. Then:
$$
\left\{ \begin{array}{llll} \e(-1) = -1,~ \e(5) = 1 ~~ & \mbox{and} ~~ & \left\{ \begin{array}{ll} n_2 = 5 & \mbox{if} ~~ c_2 \le 2 \\
n_2 = 2c_2 & \mbox{if} ~~ c_2\ge 3 \end{array} \right\} & \mbox{if} ~~ K_{\p_2} = M_1 ~, \\
&&& \\
\e(-1) = -1,~ \e(5) = 1 ~~ & \mbox{and} ~~ & \left\{ \begin{array}{ll} n_2 = 7 & \mbox{if} ~~ c_2\le 3 \\ n_2 = 2c_2 & \mbox{if}
~~ c_2\ge 4 \end{array} \right\} & \mbox{if} ~~ K_{\p_2} = M_2 ~, \\
&&& \\
\e = 1 ~~ & \mbox{and} ~~ & \left\{ \begin{array}{ll} n_2 = 7 & \mbox{if} ~~ c_2\le 3 \\ n_2 = 2c_2 & \mbox{if} ~~
c_2\ge 4 \end{array} \right\} & \mbox{if} ~~ K_{\p_2} = M_3 ~, \\
&&& \\
\e = 1 ~~ & \mbox{and} ~~ & \left\{ \begin{array}{ll} n_2 = 3 & \mbox{if} ~~ c_2\le 1 \\ n_2 = 2c_2 & \mbox{if} ~~ c_2\ge 2
\end{array} \right\} & \mbox{if} ~~ K_{\p_2} = M_4 ~.
\end{array} \right.
$$

\subsection{Degeneration of local conductors.} \label{deg} We now determine how the conductor of $r_p \otimes \phi_p$
degenerates upon reduction mod $\lambda$, i.e., we determine the number
$$
\delta_p(\phi_p) = \delta_p(K_{\p}/{\Q}_p,\phi_p) := c_p(r_p \otimes \phi_p) - c_p(\bar r_p \otimes \bar\phi_p)
$$
for characters $\phi_p \colon U_p \rightarrow \C^{\times}$ of order prime to $\ell$. Recall that $\bar r$ is the mod $\lambda$
reduction of $r$, and that $\bar r_p$ is the restriction of $\bar r$ to a decomposition group at $p$. Thus, $\bar r_p$ is the mod
$\lambda$ reduction of $r_p$.

There is a well-developed understanding of how local conductors may degenerate, cf. \cite{liv}, \cite{car}, but in our special
situation it will be more convenient and in fact basically trivial to determine $\delta_p(\phi_p)$ directly by `inspection'.
Recall that if $p$ and $\ell$ are distinct primes, and if $R\colon G_{\Q_p} \rightarrow \GL_n(\overline\Q_{\ell})$ is a
continuous $\lambda$-adic representation with $\bar{R} \colon G_{\Q_p} \rightarrow \GL_n(\overline\F_{\ell})$ a mod $\lambda$
reduction, then
$$
c_p(R) - c_p(\bar{R}) = \dim \bar{R}^{I_p} - \dim R^{I_p}
$$
where $\bar{R}^{I_p}$ and $R^{I_p}$ denote the space of fixed points of $\bar{R}(I_p)$ and $R(I_p)$ respectively (in the
associated underlying vector space). Confer \cite{liv}, section $1$.
\smallskip

We split up into 3 main subcases corresponding to the behavior of $\bar r_p$ on the inertia group $I_p$.
\smallskip

\noindent (1) Suppose that $\bar r_p$ is irreducible on $I_p$. This case will occur if $p=2$ and $r_2$ is of $A_4$- or
$S_4$-type, or if $\Ima \pi_p = \Ima \Proj(r_p)$ is totally ramified dihedral of order $2m$ with $\ell \nmid m$.

Then $\bar r_p \otimes \bar\phi_p$ and $r_p \otimes \phi_p$ are also both irreducible on $I_p$. Thus,
$$
\delta_p(\phi_p) = 0
$$
in these cases.
\smallskip

\noindent (2) Suppose that $\bar r_p$ is reducible but semisimple on $I_p$. This case occurs if $\Ima \pi_p = \Ima \Proj(r_p)$ is
either cyclic of order not divisible by $\ell$, or is ramified, but not totally ramified, dihedral of order $2m$ with $m$ not
divisible by $\ell$.

(2a) In the first case we have
$$
(r_p)_{\mid I_p} \sim \left( \begin{smallmatrix} \alpha & 0 \\ 0 & 1 \end{smallmatrix} \right) ,\quad (\bar r_p)_{\mid I_p} \sim
\left( \begin{smallmatrix} \bar\alpha & 0 \\ 0 & 1 \end{smallmatrix} \right),
$$
with $\alpha$ a ramified character of order not divisible by $\ell$. We obtain
$$
\delta_p(\phi_p) = 0
$$
since \dash--- as $\alpha$ and $\phi_p$ have order prime to $\ell$ \dash--- the characters $\alpha \phi_p$ and $\bar\alpha
\bar\phi_p$ are seen to have the same conductors, and similarly for $\phi_p$ and $\bar\phi_p$.

(2b) In the second case we have $r_p = \Ind_{M/\Q_p} (\beta)$ where $M/\Q_p$ is the unramified quadratic extension, and $\beta$
is a ramified character on $G_M$ that may be \dash--- and in fact has been \dash--- chosen to have order prime to $\ell$. Then if
$\sigma$ denotes the non-trivial automorphism of $M/\Q_p$ we have
$$
(r_p \otimes \phi_p)_{\mid I_p} \sim \left( \begin{smallmatrix} \beta \cdot (\phi_p)_{\mid G_M} & 0 \\ 0 & \beta^{\sigma} \cdot
(\phi_p)_{\mid G_M} \end{smallmatrix} \right) ,
$$
$$
(\bar r_p \otimes \bar\phi_p)_{\mid I_p} \sim \left( \begin{smallmatrix} \bar\beta \cdot (\bar\phi_p)_{\mid G_M} & 0 \\
0 & \bar\beta^{\sigma} \cdot (\bar\phi_p)_{\mid G_M} \end{smallmatrix} \right)
$$
and again we deduce
$$
\delta_p(\phi_p) = 0 ~.
$$
\smallskip

\noindent (3) Suppose finally that $\bar r_p$ is reducible but not semisimple on $I_p$. Then $\Ima \pi_p = \Ima \Proj(r_p)$ is
either cyclic of order $m$ divisible by $\ell$, or is dihedral of order $2m$ with $m$ divisible by $\ell$. Since $\Ima \pi_p$ is
isomorphic to a subgroup of one of the groups $A_4$, $S_4$, or $A_5$, these cases can only occur if $\ell \in \{ 3,5\}$ and
$m=\ell$.

(3a) In the first case we have chosen $r_p$ such that $\Ima r_p$ has the same order as $\Ima \pi_p$, i.e.,
$$
r_p \sim \left( \begin{smallmatrix} \alpha & 0 \\ 0 & 1 \end{smallmatrix} \right)
$$
with $\alpha$ a character of order $\ell$. Hence a priori
$$
\dim (r_p \otimes \phi_p)^{I_p} = \left\{ \begin{array}{ll} 1 & \mbox{ if } \phi_p \in \{ 1,\alpha^{-1} \} \\
0 & \mbox{ otherwise,} \end{array} \right.
$$
where however the case $\phi_p = \alpha^{-1}$ actually does not occur because $\phi_p$ has order prime to $\ell$ whereas $\alpha$
has order $\ell$.

On the other hand,
$$
\bar r_p \sim \left( \begin{smallmatrix} 1 & \ast \\ 0 & 1 \end{smallmatrix} \right)
$$
with $\ast \not= 0$ and so
$$
\dim (\bar r_p \otimes \bar\phi_p)^{I_p} = \left\{ \begin{array}{ll} 1 & \mbox{ if } \bar\phi_p =1 \\ 0 & \mbox{ otherwise.}
\end{array} \right.
$$
We obtain
$$
\delta_p(\phi_p) =0 ~.
$$

(3b) In the second case we have as in (2) above that $r_p = \Ind_{M/\Q_p} (\beta)$ with $M/\Q_p$ the unramified quadratic
extension and $\beta$ a ramified character on $G_M$. Hence,
$$
(r_p \otimes \phi_p)_{\mid I_p} \sim \left( \begin{smallmatrix} \beta \cdot (\phi_p)_{\mid G_M} & 0 \\ 0 & \beta^{\sigma} \cdot
(\phi_p)_{\mid G_M} \end{smallmatrix} \right)
$$
with $\sigma$ the non-trivial automorphism of $M/\Q_p$. Noticing that $(\phi_p)_{\mid G_M}$ is fixed under conjugation with
$\sigma$ whereas $\beta$ is not ($\beta^{-1}\beta^{\sigma}$ is a character of order $m=\ell$), we obtain
$$
\dim (r_p \otimes \phi_p)^{I_p} = 0 ~.
$$
On the other hand, denoting by $\alpha$ the restriction to $I_p$ of the quadratic character corresponding to $M/\Q_p$, we have
$$
(\bar r_p \otimes \bar\phi_p)_{\mid I_p} \sim \left( \begin{smallmatrix} \bar\alpha \bar\phi_p & \ast \\ 0 & \bar\phi_p
\end{smallmatrix} \right)
$$
so that
$$
\dim (\bar r_p \otimes \bar\phi_p)^{I_p} = \left\{ \begin{array}{ll} 1 & \mbox{ if } \bar\phi_p = \bar\alpha \\
0 & \mbox{ otherwise.} \end{array} \right.
$$
Consequently,
$$
\delta_p(\phi_p) = \left\{ \begin{array}{ll} 1 & \mbox{ if } \bar\phi_p = \bar\alpha \\
0 & \mbox{ otherwise.} \end{array} \right.
$$

\subsection{Local lifts: The weight.} \label{weight} In this section we define in dependence on the extension
$K_{\mathfrak l}/\Q_{\ell}$ a set $\W_{\ell}$ of natural numbers so that $k\in \W_{\ell}$ may be deduced. Recall that
$K_{\mathfrak l}/\Q_{\ell}$ is the local extension cut out by the restriction to the decomposition group $G_{\Q_{\ell}} \le
G_{\Q}$ of $\Proj (\rho)$ with $\rho$ our given modular mod $\lambda$ representation.
\smallskip

We start by quoting two theorems about the structure of the restriction $\rho_{\ell}$ of $\rho$ to the decomposition group
$G_{\Q_{\ell}} \le G_{\Q}$. Let the corresponding mod $\lambda$ cusp form be denoted by $f$. Thus $f$ is a form of weight $k\in
\{ 2,\ldots ,\ell-1 \}$. Denote by $a_n\in \overline\F_{\ell}$ the Fourier coefficients of $f$. The first theorem is due to
Deligne; a proof can be found in \cite{gro}. The statement is that if $a_{\ell} \not= 0$ then $\rho_{\ell}$ is reducible and more
precisely has the shape
$$
\rho_{\ell} \sim \left( \begin{matrix} u\cdot \chi^{k-1} & \ast \\ 0 & v \end{matrix} \right)
$$
where $u$ and $v$ are unramified characters. The second theorem \dash--- due to Fontaine with a proof to be found in \cite{edi}
\dash--- states that if $a_{\ell} =0$ then $\rho_{\ell}$ is irreducible and
$$
(\rho_{\ell})_{\mid I_{\ell}} \sim \left( \begin{matrix} (\psi^{\prime})^{k-1} & 0 \\ 0 & \psi^{k-1} \end{matrix} \right)
$$
where $\psi, \psi^{\prime} \colon I_{\ell} \rightarrow \overline\F_{\ell}$ are the 2 fundamental characters of level 2. We have
$\psi^{\prime} = \psi^{\ell}$.

As usual, we refer to these two cases as the ordinary and the supersingular case respectively.
\smallskip

Now on the other hand we have defined above in section \ref{local} a certain representation $r \colon G_{\Q} \longrightarrow
{\GL}_2(\C)$ such that in particular:
$$
\rho_{\ell} \cong \bar r_{\ell} \otimes \chi^i \otimes (\mbox{unramified character}),
$$
for some $i$, where $\bar r_{\ell}$ denotes the restriction of the mod $\lambda$ reduction of $r$ to the decomposition group
$G_{\Q_{\ell}} \le G_{\Q}$. We proceed now to compare these two descriptions of $\rho$ in order to make statements about $k$.
\smallskip

The first observation is that if $K_{\mathfrak l}/\Q_{\ell}$ were unramified, we would obviously be in the ordinary case whence
the assumption $2\le k\le \ell-1$ would imply the contradiction $k=1$. So we may start by defining
$$
\W_{\ell} :=\emptyset \quad \mbox{if } K_{\mathfrak l}/\Q_{\ell} \mbox{ unramified}.
$$
Secondly, as $\ell$ is odd the extension $K_{\mathfrak l}/\Q_{\ell}$ is either cyclic or dihedral. We split up the discussion
into 5 subcases.

\subsubsection{Tamely ramified cyclic case.} Suppose that $K_{\mathfrak l}/\Q_{\ell}$ is tamely ramified and cyclic of order
$m\in \{ 2,3,4,5 \}$. Let $e \mid m$ be the ramification index of $K_{\mathfrak l}/\Q_{\ell}$. Thus, $\ell \equiv 1 \mod e$. The
representation $\bar r_{\ell}$ has image of order $m$, and as $\ell \nmid m$ we have that $\bar r_{\ell}$ is semisimple and
reducible. We are thus in the ordinary case, and so $\chi^{k-1}$ necessarily has order $e$. On the other hand, this order is
$\frac{\ell-1}{(\ell-1,k-1)}$ so we see that
$$
k-1 = a\cdot \frac{\ell-1}{e}
$$
for some natural number $a$ which must be less than $e$, as $k\le \ell-1$.
\smallskip

So, if $e=2,3,5$ we have $a=1$, $a\in \{ 1,2\}$ and $a\in \{ 1,2,3,4\}$ respectively. If $e=4$ we must have $a\in \{ 1,3\}$ as
$a=2$ would imply $O(\chi^{k-1})=2$. So we may define
$$
\W_{\ell} = \W_{\ell}(K_{\mathfrak l}/\Q_{\ell}) := \left\{ \begin{array}{ll} \{ \frac{\ell+1}{2} \} & \mbox{  for } e=2 \\
& \\
\{ \frac{\ell+2}{3} ,\frac{2\ell +1}{3} \} & \mbox{  for } e=3 \\
& \\
\{ \frac{\ell+3}{4} ,\frac{3\ell+1}{4} \} & \mbox{  for } e=4 \\
& \\
\{ \frac{\ell+4}{5}, \frac{2\ell+3}{5}, \frac{3\ell+2}{5}, \frac{4\ell+1}{5} \} & \mbox{  for } e=5.
\end{array} \right.
$$

\subsubsection{Wildly ramified cyclic case.} If $K_{\mathfrak l}/\Q_{\ell}$ is wildly ramified cyclic of order $m$, we have
that $\bar r_{\ell}$ is twist-equivalent to
$$
\left( \begin{matrix} u & \ast \\ 0 & v \end{matrix} \right)
$$
with unramified characters $u,v$ and $\ast \not= 0$. We deduce the contradiction $k=1$ and may thus define $\W_{\ell}
:=\emptyset$ in this case.

\subsubsection{Wildly ramified dihedral cases.}\label{weight_wild_ram_dih} If $K_{\mathfrak l}/\Q_{\ell}$ is wildly ramified
dihedral of order $2m$ then, as $\ell$ is odd, we must have $\ell \mid m$. As the Galois group of $K_{\mathfrak l}/\Q_{\ell}$ is
isomorphic to a subgroup of one of $A_4$, $S_4$, $A_5$, we conclude $\ell \in \{ 3,5\}$ and $m=\ell$.

Now, $K_{\mathfrak l}/\Q_{\ell}$ has a unique quadratic subextension $M/\Q_{\ell}$, and the representation $r \colon G_{\Q}
\longrightarrow {\GL}_2(\C)$ has been chosen in such a way that $\bar r_{\ell}$ appears in shape
$$
\bar r_{\ell} \sim \left( \begin{matrix} \alpha & \ast \\ 0 & 1 \end{matrix} \right) \otimes (\mbox{unramified character})
$$
where $\alpha$ is the quadratic character corresponding to $M/\Q_{\ell}$.

So, if $M/\Q_{\ell}$ is unramified, that is, if $K_{\mathfrak l}/\Q_{\ell}$ is wildly but not totally ramified, we deduce the
contradiction $k=1$.

On the other hand, if $M/\Q_{\ell}$ is ramified, i.e., if $K_{\mathfrak l}/\Q_{\ell}$ is totally and wildly ramified, we must
have $\alpha = \chi^{k-1}$. As $k\le \ell-1$ and $\alpha$ has order 2, we deduce $k=\frac{\ell+1}{2}$. For the proof of part (b)
of Theorem \ref{t1} in the next subsection, it will be convenient at this point to make an additional remark in the special case
that $\ell=3$: We have then $k=2$ and this must be the weight attached to $\rho$ by Serre-Edixhoven; referring specifically to
Definition 4.3 of \cite{edi}, we are in case 2(b) of that definition and $k=2$ means that $K_{\mathfrak l}/\Q_{\ell}$ is `peu
ramifi\'{e}' in the sense of Serre \cite{ser1}, p.186. We may thus define
$$
\W_{\ell} := \left\{ \begin{array}{ll} \emptyset & \mbox{ for } K_{\mathfrak l}/\Q_{\ell} \mbox{ wildly but not totally ramified} \\
\{ \frac{\ell+1}{2} \} & \mbox{ if } \ell =5 \mbox{ and } K_{\mathfrak l}/\Q_{\ell} \mbox{ totally and wildly ramified,}
\end{array} \right.
$$
and for the case $\ell=3$
$$
\W_3 :=\{2 \} \quad \mbox{for } K_{\mathfrak l}/\Q_3 \mbox{ totally, wildly ramified, but peu ramifi\'{e}.}
$$

For the proof of part (b) of Theorem \ref{t1} it will then be convenient to define additionally for the special case
$\ell=3$:
$$
\W_3 := \{ 4\} \quad \mbox{for } K_{\mathfrak l}/\Q_3 \mbox{ totally, wildly ramified, and tr\`{e}s ramifi\'{e}.}
$$

\subsubsection{Tamely ramified dihedral case.}\label{tam_ram_dih} Suppose finally that $K_{\mathfrak l}/\Q_{\ell}$ is tamely ramified dihedral of order
$2m$, $m\in \{ 2,3,4,5\}$. Then $K_{\mathfrak l}$ is cyclic of degree $m$ over the unramified quadratic extension $M/\Q_{\ell}$
and $\ell \equiv -1 \mod m$. The representation $\bar r_{\ell}$ has the form $\bar r_{\ell} = \Ind_{M/\Q_{\ell}} (\beta)$ for a
tamely ramified character $\beta$ on $G_M$. Thus,
$$
(\bar r_{\ell})_{\mid I_{\ell}} \sim \left( \begin{matrix} \beta & 0 \\ 0 & \beta^{\sigma} \end{matrix} \right)
$$
where $\sigma$ denotes the non-trivial automorphism of $M/\Q_{\ell}$. We are then necessarily in the supersingular case. The
characters $\beta,\beta^{\sigma}$ are of level 2, we have
$$
\beta^{\sigma} = \beta^{\ell}
$$
and $m$ is the order of the character $\beta^{\sigma} \beta^{-1} = \beta^{\ell-1}$.

Denote as above by $\psi, \psi^{\prime}$ the two fundamental characters of level 2. Interchanging $\psi$ and $\psi^{\prime}$ if
necessary we may write
$$
\beta = \psi^r (\psi^{\prime})^s,\quad \beta^{\sigma} = \psi^s (\psi^{\prime})^r
$$
for integers $r,s$ with $0\le r < s\le \ell-1$. Then
$$
(\bar r_{\ell} \otimes \chi^{-r})_{\mid I_{\ell}} \sim \left( \begin{matrix} (\psi^{\prime})^a & 0 \\
0 & \psi^a \end{matrix} \right)
$$
with $a:=s-r$ which must satisfy $1\le a\le \ell-2$ because $\bar r_{\ell}$ is irreducible. Now,
$$
(\psi^{\prime})^a \psi^{-a}= \psi^{a\cdot (\ell-1)}
$$
must still have order $m$ and this allows us to determine the possibilities for $a$: As $\psi$ has order $\ell^2-1$ the number
$a$ must have the form
$$
a=a_0\cdot \frac{\ell+1}{m}
$$
with $1\le a_0<m$. If $m=4$ then $a_0\not= 2$ because $\psi^{a\cdot (\ell-1)}$ would otherwise have order 2.
\smallskip

On the other hand we know that the restriction $(\rho_{\ell})_{\mid I_{\ell}}$ has the form
$$
(\rho_{\ell})_{\mid I_{\ell}} \sim \left( \begin{matrix} (\psi^{\prime})^{k-1} & 0 \\
0 & \psi^{k-1} \end{matrix} \right)
$$
and is the twist by a power of $\chi =\psi^{\ell+1}$ of the above $(\bar r_{\ell})_{\mid I_{\ell}}$. We deduce that $k-1$ has the
form $a+i(\ell+1)$ or $\ell a +i(\ell+1)$ for some $i$, and since $k-1$ is a number in $\{ 1,\ldots ,\ell-2 \}$ we find that
$k=a+1$ or $k=\ell+2-a$ where $a$ is one of the numbers given above. We can then check that $k\in \W_{\ell}$ if we define the set
$\W_{\ell}$ thus:
$$
\W_{\ell} = \W_{\ell}(K_{\mathfrak l}/\Q_{\ell}) := \left\{ \begin{array}{ll}
\{ \frac{\ell+3}{2} \} & \mbox{  for  } m=2 \\ & \\ \{
\frac{\ell+4}{3} , \frac{2\ell+5}{3} \} & \mbox{  for  } m=3 \\ & \\
\{ \frac{\ell+5}{4} , \frac{3\ell+7}{4} \} & \mbox{  for  } m=4 \\ & \\
\{ \frac{\ell+6}{5} , \frac{2\ell+7}{5} , \frac{3\ell+8}{5} ,
\frac{4\ell+9}{5} \} & \mbox{  for  } m=5. \end{array} \right.
$$

\subsection{} This finishes the proof of part (a) of Theorem \ref{t1}: Starting with the modular mod $\lambda$ representation
$\rho$ cutting out $K/\Q$ with Serre invariants $(N,k,\nu)$, we defined a global character $\phi$ of order prime to $\ell$, and
then in sections \ref{local} -- \ref{weight} defined the numbers $n_p(K_{\p}/{\Q}_p , \phi_p)$ and
$\delta_p(K_{\p}/{\Q}_p,\phi_p)$, the characters $\e_p(K_{\p}/{\Q}_p)$, and the set $\W _{\ell}$, and in the process verified the
conclusions of part (a) of Theorem \ref{t1}.

\subsection{End of proof.} \label{end} We now finish the proof of Theorem \ref{t1} by proving part (b). So, let $\phi$, $N$, $k$ and
$\nu$ be given with the properties and relations stated in the theorem. In particular we have that $k$ is a number in the set
$\W_{\ell}= \W_{\ell}(K_{\mathfrak l}/\Q_{\ell})$, and because of our definition of $\W_{\ell}$ we have that $K$ is ramified over
$\ell$. Choose an embedding $G := \mathrm{Gal}(K/\Q) \hookrightarrow \PGL_2(\C)$; this gives us a projective Galois
representation
$$
\pi \colon G_{\Q} \longrightarrow \mathrm{PGL}_2(\C) ~.
$$
Denote by
$$
r \colon G_{\Q} \longrightarrow {\GL}_2(\C)
$$
the lift of $\pi$ that we specified in section \ref{local} above. The hypothesis that $K$ not be totally real implies that $r$ is
an odd representation. So, $r$ is modular of weight 1 by Langlands-Tunnell if $G\not \cong A_5$, and by assumption if $G\cong
A_5$; cf. remark (3) immediately after the statement of Theorem \ref{t1}.

Accordingly, the representation $\bar r$ defined as the reduction mod $\lambda$ of $r$ is an irreducible modular mod $\lambda$
representation. The same is then true of the representation
$$
\rho := \bar r \otimes \bar\phi
$$
with $\bar\phi$ the reduction mod $\lambda$ of $\phi$. By construction $\rho$ has conductor $N$ and character $\nu$. Now it
follows from Ribet's fundamental work \cite{rib} on `level-lowering' as complemented by Diamond in \cite{dia}, Corollary 1.2,
that $\rho$ is the mod $\lambda$ representation attached to a mod $\lambda$ cusp form of type $(N,k^{\prime},\nu)$ where
$k^{\prime}$ is the weight attached to $\rho$ by Serre as modified in \cite{edi}. Tensoring with some power of $\chi$ if
necessary, we may further assume that $k^{\prime} \le \ell+1$.
\medskip

We shall now proceed to show that $k^{\prime} \in \W_{\ell}(K_{\mathfrak l}/\Q_{\ell})$, though possibly only after tensoring
with some power of $\chi$. First notice that if $k^{\prime}=1$ then $\rho$ is unramified at $\ell$ by \cite{cv}, Corollary 0.2;
this is contrary to our assumptions, so $k^{\prime} \ge 2$.
\smallskip

Suppose that $k^{\prime} = \ell+1$. Going through the various cases of Definition 4.3 of \cite{edi}, we see that we must be in
case 2.(b) of that definition, i.e.,
$$
\rho_{\mid I_{\ell}} \sim \left( \begin{array}{cc} \chi^{\beta} & \ast \\ 0 & \chi^{\alpha} \end{array} \right) ~,
$$
where $\ast \not=0$. As we have $k^{\prime} = \ell+1$ the definition of $k^{\prime}$ is this case shows us that we must have
$a:=\min \{ \alpha ,\beta\} =0$, $b:=\max \{ \alpha ,\beta\} =1$, and that $\rho$ is not finite at $\ell$. In particular, we see
that $\Proj(\rho) (I_{\ell})$ has order $\ell \cdot (\ell-1)$. On the other hand, this group is a subgroup of one of the groups
$A_4$, $S_4$, or $A_5$, so we deduce $\ell=3$ and that this subgroup is in fact dihedral of order 6; we also see \dash--- again
by order considerations \dash--- that the projectivisation of $\rho$ is totally (wildly) ramified at $\ell$. So, $K_{\mathfrak
l}$ is a totally ramified dihedral extension of $\Q_3$ of order 6. Since $\rho$ is not finite at $\ell$, by definition this
dihedral extension is `tr\`{e}s ramifi\'{e}' in the sense of Serre, \cite{ser1}, p.186. By definition of $\W_{\ell}(K_{\mathfrak
l}/\Q_{\ell})$ in \ref{weight_wild_ram_dih} above, we have then $\W_{\ell} \ni 4 = k^{\prime}$.
\smallskip

Suppose then that $k^{\prime} = \ell$ and that we are in the ordinary case. By the theorem of Deligne quoted above in section
\ref{weight} we have then
$$
\rho_{\ell} := (\rho)_{\mid G_{\Q_{\ell}}} \sim \left( \begin{matrix} u & \ast \\ 0 & v \end{matrix} \right)
$$
for some unramified characters $u,v$. Then $\ast \not= 0$ again because $\rho$ is ramified at $\ell$. We deduce that $\ell \in \{
3,5\}$ and that $K_{\mathfrak l}/\Q_{\ell}$ is either wildly ramified cyclic of order $\ell$, or wildly but not totally ramified
dihedral of order $2\ell$. Again, this is in contradiction to our assumption $k\in \W_{\ell}$ and our definition $\W_{\ell} :=
\emptyset$ in these cases (i.e., we have explicitly excluded these cases via the assumptions of part (b) of Theorem \ref{t1}).
\smallskip

Assume then $k^{\prime} = \ell$ but that we are in the supersingular case so that $\rho_{\ell}$ has a projective kernel field
which is tamely ramified of order $2m$, $m\in \{ 2,3,4,5\}$. Now, the theorem of Fontaine quoted above implies that we have
$$
(\rho_{\ell})_{\mid I_{\ell}} \sim \left( \begin{matrix} (\psi^{\prime})^{\ell-1} & 0 \\ 0 & \psi^{\ell-1} \end{matrix} \right) ,
$$
and we deduce as above in section \ref{weight} that $\psi^{(\ell-1)^2}$ has order $m$. On the other hand, this order is
$\frac{\ell+1}{2}$, so $(m,\ell)$ is one of the following: (2,3), (3,5), (4,7). Now, it follows from \cite{edi}, Proposition 3.3,
that for a suitable power $\chi^j$ of $\chi$, the representation $\rho \otimes \chi^j$ is modular of weight 3. As one checks that
we in fact have $3\in \W_{\ell}$ in each of the $3$ cases, we are done.
\medskip

We may now conclude \dash--- possibly after tensoring $\rho$ with a suitable power of $\chi$ \dash--- that either $k^{\prime} \in
\W_{\ell}$ or else $2\le k^{\prime} \le \ell-1$. But in the latter case we also have $k^{\prime} \in \W_{\ell}$ by the already
proven part (a) of the theorem.
\medskip

We thus see that our desired conclusion, i.e., the claim of existence of a mod $\lambda$ cusp form of type $(N,k,\nu)$ with
associated projective kernel field $K$ is true for at least one number $k$ in the set $\W_{\ell}(K_{\mathfrak l}/\Q_{\ell})$. In
order to finish the proof it then clearly suffices to show the following two statements:
\smallskip

(1) In case either $\ell=5$ or $K_{\mathfrak l}/{\Q}_{\ell}$ has degree not divisible by 5, we must show that $(N,\tilde{k},\nu)$
`occurs' in the above sense for every number $\tilde{k} \in \W_{\ell}$ with $2\le \tilde{k} \le \ell-1$.
\smallskip

(2) In case $\ell \not= 5$ and $K_{\mathfrak l}/{\Q}_{\ell}$ has degree divisible by 5, we must define a partition $\W_{\ell} =
\W_{\ell}^{+} \cup \W_{\ell}^{-}$ and prove the existence of a sign $\mu \in \{ +,-\}$ such that $(N,\tilde{k},\nu)$ occurs for
every number $\tilde{k} \in \W_{\ell}^{\mu}$ with $2\le \tilde{k} \le \ell-1$.
\smallskip

This is essentially accomplished via the theory of companion forms \cite{gro}, \cite{cv}, and the theory of $\theta$-cycles
\cite{joc}, \cite{edi}, applied to the ordinary and supersingular cases respectively: Suppose that $|\W_{\ell}|>1$, and that $f$
is a mod $\lambda$ cusp form of type $(N,k,\nu)$ where $k\in \W_{\ell}$. Let $\rho_f$ denote the mod $\lambda$ representation
attached to $f$.
\smallskip

Suppose first that $K_{\mathfrak l}/\Q_{\ell}$ is tamely ramified cyclic of degree $m$ and ramification index $e\mid m$. As
$|\W_{\ell}|>1$ we have $e\in \{ 3,4,5\}$. According to the main results of \cite{gro} and \cite{cv}, the form $f$ has a
companion form of weight $k_1:=\ell+1-k$. One checks that if $e$ is 3 or 4 then $\W_{\ell}$ consists precisely of the numbers $k$
and $k_1$. If $e=5$, we define
$$
\W_{\ell}^{+} := \left\{ \frac{\ell+4}{5} , \frac{4\ell+1}{5} \right\} , \quad \W_{\ell}^{-} :=
\left\{ \frac{2\ell+3}{5} , \frac{3\ell+2}{5} \right\}
$$
so that
$$
\W_{\ell} = \W_{\ell}^{+} \cup \W_{\ell}^{-} ~,
$$
and one checks that if a subset $\W_{\ell}^{\mu}$ contains $k$, then it also contains $k_1$.
\smallskip

Suppose then that $K_{\mathfrak l}/\Q_{\ell}$ is tamely ramified dihedral of order $2m$, so that we are in the supersingular
case. As $|\W_{\ell}|>1$ we have $m\in \{ 3,4,5\}$. We then first observe that any number in $\W_{\ell}$ is $\ge 3$ except if
$\ell=3$ and $m=4$; however, in the latter case we have $\W_{\ell} = \{ 2,4\}$ and there is nothing to show, as $4> \ell-1$ in
this case. So, we may assume $k\ge 3$. Now the theory of $\theta$-cycles \cite{joc}, \cite{edi}, in particular Proposition 3.3 in
\cite{edi}, implies that $\rho_f \otimes \chi^{\ell+1-k}$ is modular of type $(N,\ell+3-k,\nu)$ (and obviously with $K$ its
projective kernel field). One checks that for $m=$ 3 or 4, the sum of the 2 numbers in $\W_{\ell}$ is indeed $\ell+3$, so we are
done in those cases. If $m=5$ we can split $\W_{\ell}$ into 2 subsets:
$$
\W_{\ell} = \W_{\ell}^{+} \cup \W_{\ell}^{-}
$$
with
$$
\W_{\ell}^{+} := \left\{ \frac{\ell+6}{5} , \frac{4\ell+9}{5} \right\} ,\quad \W_{\ell}^{-} := \left\{ \frac{2\ell+7}{5} ,
\frac{3\ell+8}{5} \right\} .
$$
As the sum of the two numbers in each $\W_{\ell}^{\mu}$ is $\ell+3$, we are done.
\medskip

Thus, the proof of part (b) of Theorem \ref{t1} is finished.

\section{Additional remarks.} \label{rem}

\subsection{} The purpose of Theorem \ref{t1} is to give the possible Serre invariants $(N,k,\nu)$ of mod $\ell$
representations of one of the exceptional types $A_4$, $S_4$, or $A_5$ in terms of local data attached to a fixed projective
kernel field $K$. One may of course reverse the question, i.e., start with a triple $(N,k,\nu)$ and ask whether there are any
cusp forms with these data and an attached exceptional mod $\ell$ representation. In connection with this question, the Theorem
will immediately give information about whether there are any local Galois theoretic obstructions for this to be the case, and if
there are not, will give detailed information about the local structure of possible corresponding projective kernel fields. This
information will suffice to determine all possible candidates for these projective kernel fields, either via class field theory
in the $A_4$- and $S_4$-cases, or via geometry of numbers (or tables) in the $A_5$-case. This gives an algorithm which will
decide by a finite amount of computation the above question, \dash--- again of course under the assumption of the Artin
conjecture in the $A_5$ case. Examples of this were given above in section \ref{exam}. In these examples the dimension of the
space of cusp forms in question was sufficiently low that we could actually `point' to the exceptional form. In general however,
the theorem will only inform us about the existence of an exceptional mod $\ell$ cusp form with the given data $(N,k,\nu)$,
though \dash--- as was pointed out above \dash--- we may in fact count the number of such forms. If one wishes additionally to
`point' to these forms one must first realize that all this can mean in general is to give an algorithm that computes their
Fourier expansions. We shall now briefly indicate how this can be done.
\smallskip

Starting with an extension $K/\Q$ of $A_4$, $S_4$ or $A_5$ type, we choose an embedding $\Gal(K/\Q) \hookrightarrow \PGL_2(\C)$
and obtain thus a complex projective Galois representation $\pi$. The Fourier expansion that we desire is obtained via reduction
mod $\lambda$ \dash--- for some prime $\lambda$ over $\ell$ \dash--- from the L-series of a twist of a lift of $\pi$ to a linear
representation $G_{\Q} \longrightarrow \GL_2(\C)$. If one uses the results the papers \cite{que}, \cite{cre}, \cite{cre1}, and
\cite{cre2}, one has an algorithm that computes the L-series of {\it some} lift $r_1$ of $\pi$ (corresponding to a central
embedding problem with cyclic kernel of 2-power order). The desired Fourier expansions can be computed from twists of the mod
$\lambda$ reduction $\rho_1$ of $r_1$. Now the problem is that one does not a priori know how the conductor of $\rho_1$ behaves
under twists. However, this can easily be determined by Theorem \ref{t1} or rather the proof if it: The representation $\rho_1$
is a twist of the mod $\lambda$ reduction $\rho$ of the special linear lift $r$ occurring in the proof of Theorem \ref{t1}
(section \ref{main} above), say
$$
\rho_1 \cong \rho \otimes \xi ~.
$$
If we can determine $\xi$ then because we have complete information about the behavior of $\rho$ with respect to twists, we have
determined the Serre invariants of all twists of $\rho_1$ (essentially parametrized by the conductor of the character by which we
are twisting). But $\xi$ is easily determined: Let $\Sigma$ denote the product of $\ell$ and the finite primes at which $\rho$ or
$\rho_1$ ramifies. As $\rho_1$ is given via an explicitly solved embedding problem, one can determine $(\rho_1)_{\mid I_p}$ for
each of the finitely many primes $p\mid \Sigma$. On the other hand $(\rho)_{\mid I_p}$ is known by construction. Hence, the
quantity $\xi_p := \xi_{\mid I_p}$ can be determined for $p\mid \Sigma$. Identifying $\xi_p$ with a character on $U_p$, we have
$\xi$ explicitly determined as a mod $\ell$ Dirichlet character by
$$
\xi \colon r\mapsto \prod_{p\mid \Sigma} \xi_p(r)^{-1} ,\quad \mbox{for } (r,\Sigma) = 1 ~.
$$

We can also use the global information obtained by solving an embedding problem to resolve the indeterminacy in part (b) of
Theorem \ref{t1}: In case that $\ell \not= 5$ and $K_{\mathfrak l}/{\Q}_{\ell}$ has degree divisible by 5, the theorem informs us
\dash--- under the appropriate modularity assumption \dash--- that $\rho$ is mod $\ell$ modular of type $(N,k,\nu)$ with $k$ any
of the 2 weights in $\W_{\ell}^{\mu}$ for {\it some} choice of sign $\mu \in \{ +,-\}$. However, since we now know
$$
(\rho)_{\mid I_{\ell}} \cong (\rho_1 \otimes \xi^{-1})_{\mid I_{\ell}}
$$
explicitly, we can \dash--- as is easily checked \dash--- determine $\mu$ by looking at this restriction.

\subsection{} Continuing the above analysis, if we have given a mod $\ell$ cusp form $f$ of minimal level and weight $k\in
\{ 2,\ldots \ell-1 \}$, and if we know that the Galois representation $\rho_f$ attached to $f$ is irreducible, then we have
\dash--- under the assumption of the Artin conjecture \dash--- an effective procedure which will determine whether $\rho_f$ is of
one of the exceptional types $A_4$, $S_4$, $A_5$. If it is not, the corresponding projective kernel field has Galois group either
dihedral or isomorphic to $\PSL_2(\F_{\ell^m})$ or $\PGL_2(\F_{\ell^m})$ for some $m$. Since the dihedral case means that
$\rho_f$ is induced from a global 1-dimensional character, it is easy to see \dash--- as also noted in the introduction \dash---
that it can be effectively decided whether this case prevails. If it does not, we are in the $\PSL_2$ or $\PGL_2$ case and the
question arises how we can effectively distinguish between these, and in particular how to determine $m$ computationally. In
general, these questions lead to certain not completely trivial considerations that we shall be addressing elsewhere.

\subsection{Acknowledgements.} The authors wish to thank G. B\"{o}ckle for drawing their attention to the possibility of
using the theorem of Fong--Swan in the proof of Proposition \ref{p1} (originally we split up the discussion according to whether
$\ell \nmid \# G$ or not). This led to a nice streamlining of the argument.

\end{document}